\documentclass[a4paper,12pt]{amsart}
\usepackage{amsfonts}
\usepackage{amssymb}
\usepackage{ifthen}
\usepackage{graphicx}
\nonstopmode \numberwithin{equation}{section}
\newtheorem{thm}{Theorem}
\newtheorem{cor}{Corollary}
\newtheorem{lem}{Lemma}
\newtheorem{prop}{Proposition}
\newtheorem{opbl}{Open Problem}

\newtheorem{claim}{Claim}
\newtheorem{conj}[equation]{Conjecture}

\theoremstyle{definition}
\newtheorem{defn}{Definition}
\newtheorem{case}{Case}
\newtheorem{examp}[equation]{Example}
\newtheorem{prob}[equation]{Problem}
\newtheorem{ques}[equation]{Question}
\newtheorem{rem}{Remark}
\newcounter {own}
\def\theown {\thesection       .\arabic{own}}

\newenvironment{pf}[1][]{%
 \vskip 3mm
 \noindent
 \ifthenelse{\equal{#1}{}}%
  {{\slshape Proof. }}%
  {{\slshape #1.} }%
 }%
{\qed\bigskip}

\newcounter{alphabet}
\newcounter{tmp}
\newenvironment{Thm}[1][]{\refstepcounter{alphabet}%
\bigskip%
\noindent%
{\bf Theorem \Alph{alphabet}}%
\ifthenelse{\equal{#1}{}}{}{ (#1)}%
{\bf .} \itshape}{\vskip 8pt}

\makeatletter   % 2011-10-28 Correction
\newcommand{\Ref}[1]{\@ifundefined{r@#1}{}{\setcounter{tmp}{\ref{#1}}\Alph{tmp}}}
\makeatother

\def\be{\begin{equation}}
\def\ee{\end{equation}}

\newcommand{\bee}{\begin{enumerate}}
\newcommand{\eee}{\end{enumerate}}

\newcommand{\blem}{\begin{lem}}
\newcommand{\elem}{\end{lem}}
\newcommand{\bthm}{\begin{thm}}
\newcommand{\ethm}{\end{thm}}
\newcommand{\bcor}{\begin{cor}}
\newcommand{\ecor}{\end{cor}}
\newcommand{\beg}{\begin{examp}}
\newcommand{\eeg}{\end{examp}}
\newcommand{\begs}{\begin{examples}}
\newcommand{\eegs}{\end{examples}}
\newcommand{\bdefe}{\begin{defn}}
\newcommand{\edefe}{\end{defn}}
\newcommand{\bprob}{\begin{prob}}
\newcommand{\eprob}{\end{prob}}
\newcommand{\bques}{\begin{ques}}
\newcommand{\eques}{\end{ques}}
\newcommand{\bei}{\begin{itemize}}
\newcommand{\eei}{\end{itemize}}
\newcommand{\bop}{\begin{opbl}}
\newcommand{\eop}{\end{opbl}}
\newcommand{\bcl}{\begin{claim}}
\newcommand{\ecl}{\end{claim}}

\newcommand{\bca}{\begin{case}}
\newcommand{\eca}{\end{case}}
\newcommand{\bcon}{\begin{conj}}
\newcommand{\econ}{\end{conj}}
\newcommand{\bcons}{\begin{conjs}}
\newcommand{\econs}{\end{conjs}}
\newcommand{\bprop}{\begin{prop}}
\newcommand{\eprop}{\end{prop}}
\newcommand{\br}{\begin{rem}}
\newcommand{\er}{\end{rem}}
\newcommand{\brs}{\begin{rems}}
\newcommand{\ers}{\end{rems}}
\newcommand{\bo}{\begin{obser}}
\newcommand{\eo}{\end{obser}}
\newcommand{\bos}{\begin{obsers}}
\newcommand{\eos}{\end{obsers}}
\newcommand{\bpf}{\begin{pf}}
\newcommand{\epf}{\end{pf}}
\newcommand{\ba}{\begin{array}}
\newcommand{\ea}{\end{array}}
\newcommand{\beq}{\begin{eqnarray}}
\newcommand{\beqq}{\begin{eqnarray*}}
\newcommand{\eeq}{\end{eqnarray}}
\newcommand{\eeqq}{\end{eqnarray*}}

\newcommand{\ds}{\displaystyle}

\newcounter{minutes}
\divide\time by 60
\newcounter{hours}
\multiply\time by 60 \addtocounter{minutes}{-\time}
\begin{document}

\bibliographystyle{amsplain}

\title[Neargeodesics in John domains in Banach spaces] { Neargeodesics in John domains in Banach spaces}

%=========================================================================
\thanks{$^\dagger$ File:~\jobname .tex,
          printed: \number\year-\number\month-\number\day,
          \thehours.\ifnum\theminutes<10{0}\fi\theminutes}
%=========================================================================

\author{Yaxiang Li}
\address{Yaxiang Li,  College of Science,
Central South University of
Forestry and Technology, Changsha,  Hunan 410004, People's Republic
of China} \email{yaxiangli@163.com}

\date{}

\begin{abstract}Let $E$ be a real Banach space with dimension at least
$2$. In this paper, we prove that  if
$D\subset E$ is a John domain which is homeomorphic to an inner
uniform domain via a CQH map, then
 each neargeodesic in $D$ is a cone arc.
\end{abstract}

%\thanks{The research was partly supported by NSFs of
%China (No. 11101138 and 11071063) and Hunan Provincial Innovation Foundation For
%Postgraduate. }

\maketitle \pagestyle{myheadings} \markboth{}{Neargeodesics in John
domains in Banach spaces}

\section{Introduction}

Conformally invariant metrics, such as the hyperbolic metric, are some of the key tools of classical
function theory of plane. Quasiconformal and quasiregular mapping \cite{Vai1,Vu} generalize this theory to the Euclidean $n$-dimensional spaces. In the higher dimensions $n\geq 3$ there is no counterpart of the hyperbolic metric for a general subdomain of $\mathbb{R}^n$. In this case one can, however, introduce new metrics, hyperbolic type metrics, which still have some properties of the hyperbolic metrics. Also, it is useful to study domains where various hyperbolic type metrics compared. These ideas were presented for the first time in book form in \cite{Vu}. Recently, many authors have studied this topic \cite{HIMPS,Klen,krt,rt,rt2}.

Hyperbolic type metrics generalize also to Banach spaces. In this case for instance quasiconformality is defined in terms of the quasihyperbolic metric introduced in \cite{GP}. The present paper deals with the hyperbolic type geometries in so called John domains of Banach spaces. For the statement of our main result we introduce some terminology and notation.

Throughout the paper, we always assume that $E$ denotes a real
Banach space with dimension at least $2$. The norm of a vector $z$
in $E$ is written as $|z|$, and for each pair of points $z_1$, $z_2$
in $E$, the distance between them is denoted by $|z_1-z_2|$, the
closed line segment with endpoints $z_1$ and $z_2$ by $[z_1, z_2]$.
We always use $\mathbb{B}(x_0,r)$ to denote the open ball $\{x\in
E:\,|x-x_0|<r\}$ centered at $x_0$ with radius $r>0$. Similarly, for
the closed balls and spheres, we use the usual notations
$\overline{\mathbb{B}}(x_0,r)$ and $ \mathbb{S}(x_0,r)$,
respectively.

\bdefe\label{def-1} A domain $D$ in $E$ is called $c$-{\it John domain}
in the norm metric provided there exists a constant $c$ with the
property that each pair of points $z_{1},z_{2}$ in $D$ can be joined
by a rectifiable arc $\alpha$ in $ D$ such that for all $z\in \alpha$ the following holds:
\begin{equation} \label{eq-1} \ds\min\{\ell (\alpha [z_1, z]), \; \ell (\alpha
[z_2, z])\}\leq c\,d_{D}(z),\end{equation}
\noindent where $d_D(z)$ denotes the
distance from $z$ to the boundary $\partial D$ of $D$, $\ell(\alpha)$ denotes the length of $\alpha$,
$\alpha[z_{j},z]$ the part of $\alpha$ between $z_{j}$ and $z$ (cf.
\cite{Bro,MS, NV,Jo}). The arc $\alpha$ is called to be a {\it $c$-cone arc} . \edefe

%A domain $D$ in $E$  is said to be a {\it $c$-uniform
%domain} (cf.\cite{Martio-80,MS, Vai, Vai4}) if there
%is a constant $c\geq 1$ such that each pair of points $z_1,z_2\in D$
%can be joined by an arc $\alpha$  satisfying \eqref{eq-1} and
%
%\beq\label{eq-2} \ell(\alpha)\leq c\,|z_{1}-z_{2}|. \eeq The we say that $\alpha$ is a {\it  $c$-uniform
%arc} (cf.
%\cite{Vai4}).

For $z_1$, $z_2\in D$,  the {\it inner length metric $\lambda_D(z_1,
z_2)$} between them is defined by
$$\lambda_D(z_1,z_2)=\inf \{\ell(\alpha):\; \alpha\subset D\;
\mbox{is a rectifiable arc joining}\; z_1\; \mbox{and}\; z_2 \}.$$

\bdefe\label{def-2}A domain $D$ in $E$ is called {\it an inner $c$-uniform domain} if there
is a constant $c\geq 1$ such that each pair of points $z_1,z_2\in D$
can be joined by an arc $\alpha$  satisfying \eqref{eq-1} and
$$ \ell(\alpha)\leq c \lambda_D(z_1, z_2).$$
Such an arc
 $\alpha$ is called to be an {\it inner $c$-uniform arc} (cf. \cite{Vai7}).\edefe

From the Definition \ref{def-2}, we see that an inner $c$-uniform domain is  $c$-John. If $C$ is any compact subset of the line segment $[0,e_1]\subset\mathbb{R}^2$, then $\mathbb{B}^2(0,1)\setminus C$ is $c$-John with a universal $c$, but it need not be inner uniform. This example is due to J. Heinonen and presented by V\"ais\"al\"a in \cite{Vai7}. See \cite{BHK,HW,HLSW,Vai7} for more details on John domains and  inner uniform domains.

  In 1989, Gehring, Hag and Martio \cite{Gm} discussed the following question.
 \bques\label{ques}Suppose that $D\subset \mathbb{R}^n$ is a $c$-John domain and that $\gamma$ is a quasihyperbolic geodesic in $D$. Is $\gamma$ a $b$-cone arc for some $b=b(c)$?\eques

And they proved the following result.

 \begin{Thm}\cite[Theorem 4.1]{Gm}  If $D\subset \mathbb{R}^2$ is a simply connected John domain, then every quasihyperbolic or hyperbolic geodesic in $D$ is a cone arc. \end{Thm}

 Meanwhile, they construct several examples to show that a quasihyperbolic geodesic in a $c$-John domain need not be a $b$-cone arc with $b=b(c)$ unless $n=2$ and $D$ is simply connected.

 In 1989, Heinonen \cite{H} proposed the following question.

 \bques\label{ques1}\cite{H} Suppose that $D\subset \mathbb{R}^n$ is a $c$-John domain which  quasiconformally equivalent to the unit ball $\mathbb{B}$ and that $\gamma$ is a quasihyperbolic geodesic in $D$. Is $\gamma$ a $b$-cone arc for some constant $b$?\eques

In 2001,  Bonk, Heinonen and Koskela \cite[Theorem 7.12]{BHK} give an affirmative answer to Question \ref{ques1}. We remark that every ball is inner uniform.

 \begin{Thm}\label{ThmBHK}\cite[Proposition 7.12]{BHK}  If $D\subset \mathbb{R}^n$ is a bounded $a$-John domain which is homeomorphic to an inner $c$-uniform domain via a $K$-quasiconformal map, then each quasihyperbolic geodesic in $D$ is a $b$-cone arc with $b=b(a,c,K,n)$. \end{Thm}

We note that the constant $b$ in  Theorem \Ref{ThmBHK} depends on the dimensional $n$ and the modulus estimates of curves is used in the proof of  Theorem \Ref{ThmBHK}. As is known to all, the method of path families is useless in Banach spaces. Hence, it is natural to ask that if  Theorem \Ref{ThmBHK} could be dimensional free or not. In other words, does  it holds in Banach spaces or not. The main aim of this paper is to consider this problem. Our result shows that the answer to the problem is affirmative, and the condition ``bounded" in Theorem \Ref{ThmBHK} is redundant. our main result is as follows.

\begin{thm}\label{thm1}
Suppose that $D\subset E$ is an $a$-John domain which is
homeomorphic to an inner $c$-uniform domain via an (M,C)-CQH. Let
$z_1$, $z_2\in D$ and $\gamma$ be a $c_0$-neargeodesic joining $z_1$
and $z_2$ in $D$. Then $\gamma$ is a $b$-cone arc, where the
positive constant $b$ depends only on $a$, $c$, $c_0$, $C$ and $M$.
\end{thm}

The organization of this paper is as follows. In section \ref{sec3}, we prove several lemmas which is critical to the proof of our main result and in section \ref{sec4}, we will prove Theorem \ref{thm1}. In section \ref{sec2},  some preliminaries are stated.

\section{Preliminaries}\label{sec2}

The {\it quasihyperbolic length} of a rectifiable arc or a path
$\gamma$ in $D$ is the number (cf.
 \cite{Avv,Geo,GP,Vai3})
$$\ell_k(\gamma)=\int_{\gamma}\frac{1}{d_D(z)}\,|dz|.
$$

For each pair of points $z_1$, $z_2$ in $D$, the {\it
quasihyperbolic distance} $k_D(z_1,z_2)$ between $z_1$ and $z_2$ is
defined in the usual way:
$$k_D(z_1,z_2)=\inf\ell_k(\alpha),
$$
where the infimum is taken over all rectifiable arcs $\alpha$
joining $z_1$ to $z_2$ in $D$.

 For all $z_1$, $z_2$ in $D$, we have
(cf. \cite{Vai3})
\beq\label{eq(0000)} k_{D}(z_1, z_2)\geq
\log\Big(1+\frac{\lambda_D(z_1,z_2)}{\min\{d_{D}(z_1), d_{D}(z_2)\}}\Big)
 \geq
\Big|\log \frac{d_{D}(z_2)}{d_{D}(z_1)}\Big|.\eeq
Moreover, if $|z_1-z_2| < d_D(z_1)$, we have
\cite[Lemma 3.7]{ Vu}
\beq \label{vu1}
k_D(z_1,z_2)\leq \log\Big( 1+ \frac{
|z_1-z_2|}{d_D(z_1)-|z_1-z_2|}\Big).
\eeq

Gehring and Palka \cite{GP} introduced the quasihyperbolic metric of
a domain in $R^n$, and it has been recently used by many authors
 in the study of quasiconformal mappings and related questions \cite{Geo,HIMPS, krt, rt} etc. Recall that an arc $\alpha$ from $z_1$ to
$z_2$ is a {\it quasihyperbolic geodesic} if
$\ell_k(\alpha)=k_D(z_1,z_2)$. Obviously, each subarc of a quasihyperbolic
geodesic is a quasihyperbolic geodesic. It is known that a
quasihyperbolic geodesic between every pair of points in $E$ exists
if the dimension of $E$ is finite, see \cite[Lemma 1]{Geo}. This is
not true in infinite dimensional Banach spaces (cf. \cite[Example
2.9]{Vai3}). In order to remedy this shortage, V\"ais\"al\"a
introduced the following concept \cite{Vai4}.

\bdefe \label{def1.4}Let $D\neq E$ and $c\geq 1$. An arc
$\alpha\subset D$ is a $c$-neargeodesic if
$\ell_k(\alpha[x,y])\leq c\;k_D(x,y)$ for all $x, y\in \alpha$.
\edefe

In \cite{Vai4}, V\"ais\"al\"a proved the following property
concerning the existence of neargeodesics in $E$.

\begin{Thm}\label{LemA} $($\cite[Theorem 3.3]{Vai4}$)$
Let $\{z_1,\, z_2\}\subset D$ and $c>1$. Then there is a
$c$-neargeodesic in $D$ joining $z_1$ and $z_2$.
\end{Thm}

Now let us recall the following characterization of inner uniform
domains, which is due to V\"{a}is\"{a}l\"{a}.

\begin{Thm}\label{ThmC} $($\cite[Theorem 2.33]{Vai7}$)$ A domain
$D\subset E$ is an inner $c$-uniform domain if and only if $k_D(x,y)\leq
c'\;
 \log\Big(1+\frac{\lambda_D(x,y)}{\ds\min\{d_D(x),d_D(y)\}}\Big)$ for all $x,y\in D$,
where the constants $c$ and $c'$ depend only on each other.\end{Thm}

%uniform domains and
%Obviously, uniform domains are inner uniform domains, but inner uniform does not imply uniform.
%

Generalizing quasiconformal, V\"{a}is\"{a}l\"{a} introduced CQH
homeomorphisms (cf.  \cite{Vai1,Vai4}).

\bdefe \label{def1.6} Suppose $f: D\to D'$ is a homeomorphism. Then $f$
is said to be {\it $C$-coarsely $M$-quasihyperbolic}, or briefly $(M,C)$-CQH,
if it satisfies
$$\frac{k_D(x,y)-C}{M}\leq k_{D'}(f(x),f(y))\leq M\;k_D(x,y)+C$$
for all $x$, $y\in D$. \edefe

\section{Properties of cone arcs} \label{sec3}
In what follows, we always assume that  $f: D\to D'$ is an $(M,C)$-CQH map, that $D$ is an $a$-John domain and that $D'$ is an inner $c$-uniform domain. Also we use $x$, $y$, $z$, $\cdots$ to
denote the points in $D$, and $x'$, $y'$, $z'$, $\cdots$ the images
of $x$, $y$, $z$, $\cdots$ in $D'$, respectively, under $f$. For
arcs $\alpha$, $\beta$, $\gamma$, $\cdots$ in $D$, we also use
$\alpha'$, $\beta'$, $\gamma'$, $\cdots$ to denote their images in
$D'$.

For $x$, $y\in D$, let $\beta$ be an arc joining $x$ and $y$ in
$D$. We are going to determine some special points on $\beta'$.

\subsection{Determination of special points on $\beta'$}\label{sub3.1}

Without loss of generality, we may assume that $d_{D'}(y')\geq d_{D'}(x')$.
Then there must exist a point $w'_0\in\beta'$ which is the first
point along the direction from $x'$ to $y'$ such that

$$d_{D'}(w'_0)=\sup\limits_{p'\in \beta'}d_{D'}(p').
$$
It is possible that $w'_0=x'$ or $y'$. Obviously, there exists a
nonnegative integer $m$ such that

$$ 2^{m}\, d_{D'}(x') \leq d_{D'}(w'_0)< 2^{m+1}\, d_{D'}(x'),$$
and $x'_0$ the first point in $\beta'[x',w'_0]$ from $x'$ to $w'_0$
with

$$d_{D'}(x'_0)=2^{m}\, d_{D'}(x').
$$

Let $x_1'=x'$. If $x'_0= x'_1$, we let $x'_2=w'_0$. It is possible
that $x'_1=x'_2$. If $x'_0\not= x'_1$, then we let $x'_2,\ldots
,x'_{m+1}\in \beta'[x',x'_0]$ be the points such that for each $i\in
\{2,\ldots,m+1\}$, $x'_i$ denotes the first point from $x'$ to
$x'_0$ with
%\be\label{hws-eq(4.3)}
$$d_{D'}(x'_i)=2^{i-1}\, d_{D'}(x'_1).
$$
Obviously, $x'_{m+1}=x'_0$. If $x'_0\not= w'_0$, then we use
$x'_{m+2}$ to denote $w'_0$.

In a similar way, let $s\geq 0$ be the integer such that

$$2^{s}\, d_{D'}(y')
\leq d_{D'}(w'_0)< 2^{s+1}\, d_{D'}(y'),
$$
and $x_{1,0}'$ the first point in $\beta'[y', x_{1,0}']$ from $y'$
to $x_{1,0}'$ with

$$d_{D'}(x_{1,0}')=2^{s}\, d_{D'}(y').
$$
Let $x_{1,1}'=y'$. If $x_{1,0}'=x_{1,1}'$, we let
$x_{1,2}'=x_{1,0}'$. It is possible that $x_{1,2}'=x_{1,1}'$. If
$x_{1,0}'\not= y'$, then we let $x_{1,2}',\ldots , x_{1,s+1}'$ be
the points in $\beta'[y',w'_0]$ such that for each  $j\in
\{2,\ldots, s+1\}$, $x_{1,j}'$ is the first point from $x_{1,1}'$ to
$w'_0$ with
%\be\label{hws-eq(4.3)}
$$d_{D'}(x_{1,j}')=2^{j-1}\, d_{D'}(x_{1,1}').
$$
Then $x_{1,s+1}'=x'_{1,0}$. If $x_{1,0}'\not= w'_0$, we let
$x_{1,s+2}'=w'_0$.

\subsection{Elementary properties} \label{sub3.2}
In the following, we assume that  for each
$s_1$, $s_2\in\beta$, \be\label{eq-3} \ell_{k}(\beta[s_1,s_2])\leq
4a^2c_0k_{D}(s_1,s_2)+4a^2c_0,\ee where $a$ and $c_0$ are the same constants as in Theorem \ref{thm1}.  Obviously, (\ref{eq-3}) is
satisfied for each $c_0$-neargeodesic.

%By Proposition \ref{prop1.3}, in the following, we assume that $f:
%D\mapsto D'$ is an $(M, C)$-CQH homeomorphism, where $(M, C)$
%depends only on $(K, n)$.

\blem \label{lem13} For each $k\in \{1,\cdots,m\}$ and $z'\in
\beta'[x'_k,x'_{k+1}]$,

\bee \item $d_{D'}(x'_{k+1})\leq a_2\;d_{D'}(z')$;

\item $\lambda_{D'}(x'_{k+1},x'_k)\leq
a_2\;d_{D'}(z')$ and

\item  $\max\{\lambda_{D'}(z',x'_k), \lambda_{D'}(x'_{k+1},z')\}\leq
a_2\;d_{D'}(z')$, \eee where
$a_2=(1+2a_1)^{4a^2c_0c'M^2+1}e^{C+4a^2c_0M+4a^2c_0CM}$,
$a_1=e^{3(C+1)(a_0+M)}$ and $a_0=2^4[c'+4a^2c_0c'M+C+4a^2c_0]^4$. Here and
in what follows, $[\cdot]$ always denotes the greatest integer
part.\elem

\bpf At first, we prove the following inequality: For any $k\in \{1,
\cdots, m\}$,
\beq\label{eq(0-02)} \lambda_{D'}(x'_{k+1},x'_k)< a_1\; d_{D'}(x'_{k+1}).
\eeq We prove this inequality by contradiction. Suppose on the contrary that
 \beq\label{eq(0-1)} \lambda_{D'}(x'_{k+1},x'_k)\geq
a_1\; d_{D'}(x'_{k+1}). \eeq Let $y'_{k,1}$, $y'_{k,2}$, $\cdots$,
$y'_{k,a_0+1}\in \beta'[x'_k, x'_{k+1}]$ be $a_0+1$ points such that
$y'_{k,1}=x'_k$, $y'_{k,a_0+1}=x'_{k+1}$ and
$\lambda_{D'}(y'_{k,i+1},y'_{k,i})\geq\frac{\lambda_{D'}(x'_k,x'_{k+1})}{a_0}$.
Then for each $i\in\{1, 2, \cdots, a_0\}$,
\begin{eqnarray*}k_{D'}(y'_{k,i},y'_{k,i+1})&\geq&
\log\Big(1+\frac{\lambda_{D'}(y'_{k,i+1},y'_{k,i})}{\min\{d_{D'}(y'_{k,i+1}),d_{D'}(y'_{k,i})\}}\Big)\\
\nonumber&\geq&
\log\Big(1+\frac{\lambda_{D'}(x'_k,x'_{k+1})}{2a_0d_{D'}(x'_k)}\Big).
\end{eqnarray*}

We see from \eqref{eq-3} and Theorem \Ref{ThmC} that
\begin{eqnarray*}a_0\,\log\Big(1+\frac{\lambda_{D'}(x'_k,x'_{k+1})}{2a_0d_{D'}(x'_k)}\Big)&\leq&
\sum_{i=1}^{a_0}k_{D'}(y'_{k,i},y'_{k,i+1})\leq
M\sum_{i=1}^{a_0}k_{D}(y_{k,i},y_{k,i+1})+a_0C\\ \nonumber&\leq& M
\ell_k(\beta[x_k,x_{k+1}])+a_0C\\ \nonumber&\leq&
4a^2c_0Mk_{D}(x_k,x_{k+1})+4a^2c_0M+a_0C\\ \nonumber&\leq&
4a^2c_0M^2k_{D'}(x'_k,x'_{k+1})+(a_0+4a^2c_0M)C+4a^2c_0M
\\ \nonumber&\leq& 4a^2c_0c'M^2\log\Big(1+\frac{\lambda_{D'}(x'_k,x'_{k+1})}{d_{D'}(x'_k)}\Big)\\ \nonumber&&+(a_0+4a^2c_0M)C+4a^2c_0M,
\end{eqnarray*}
whence
$$a_0\,\log\Big(1+\frac{\lambda_{D'}(x'_k,x'_{k+1})}{2a_0d_{D'}(x'_k)}\Big)\leq
8a^2c'M^2\;\log\Big(1+\frac{\lambda_{D'}(x'_k,x'_{k+1})}{d_{D'}(x'_k)}\Big),$$
which contradicts with (\ref{eq(0-1)}). Hence (\ref{eq(0-02)})
holds.

\medskip

We infer from (\ref{eq(0-02)}) that for each $z'\in
\beta'[x'_k,x'_{k+1}]$,
\beq\label{eq(0-2')}
\log\frac{d_{D'}(x'_{k+1})}{d_{D'}(z')}&<&k_{D'}(z',x'_{k+1})
\leq M k_D(z,x_{k+1})+C\\ \nonumber &\leq&M \ell_k(\beta[x_k, x_{k+1}])+C \\ \nonumber &\leq&
4a^2c_0M\;k_D(x_k, x_{k+1})+C+4a^2c_0M\\ \nonumber
 &\leq& 4a^2c_0M^2\;k_{D'}(x'_k,x'_{k+1})+4a^2c_0CM+C+4a^2c_0M\\ \nonumber
&\leq&
4a^2c_0c'M^2\;\log\Big(1+\frac{\lambda_{D'}(x'_k,x'_{k+1})}{d_{D'}(x'_k)}\Big)\\ \nonumber&&+C+4a^2c_0CM+4a^2c_0M\\
\nonumber
 &\leq&
4a^2c_0c'M^2\;\log(1+2a_1)+C+4a^2c_0CM+4a^2c_0M,\eeq
\noindent which implies that Lemma \ref{lem13} $(1)$ holds.

Hence, (\ref{eq(0-02)}) and (\ref{eq(0-2')}) yield that
$$\lambda_{D'}(x'_k,x'_{k+1}) \leq
(1+2a_1)^{4a^2c_0c'M^2+1}e^{4a^2c_0CM+4a^2c_0M+C}\;d_{D'}(z'),$$ whence Lemma
\ref{lem13} $(2)$ follows.
\medskip

Obviously,
\begin{eqnarray*}\log\Big(1+\frac{\lambda_{D'}(x'_{k},z')}{d_{D'}(z')}\Big)
&\leq&k_{D'}(x'_k,z')
\leq M k_{D}(x_k, z)+C\\ \nonumber &\leq&M \ell_k(\beta[x_k,x_{k+1}])+C\\ \nonumber &\leq& 4a^2c_0Mk_{D}(x_k,x_{k+1})+4a^2c_0M+C\\ \nonumber &\leq&
4a^2c_0c^2\;k_{D'}(x'_k,x'_{k+1})+C+4a^2c_0M+4a^2c_0CM\\ \nonumber &\leq&
4a^2c_0c'M^2\;\log\Big(1+\frac{\lambda_{D'}(x'_k,x'_{k+1})}{d_{D'}(x'_k)}\Big)\\ \nonumber &&+C+4a^2c_0M+4a^2c_0CM,
\end{eqnarray*} which, together with (\ref{eq(0-02)}),
yields
\beq\label{eq(0-7-1)} \lambda_{D'}(x'_k,z')\leq
(1+2a_1)^{4a^2c_0c'M^2}e^{C+4a^2c_0M+4a^2c_0CM}\;d_{D'}(z').\eeq

The similar discussion as in (\ref{eq(0-7-1)}) shows that
\beq\label{eq(0-7-1')} \lambda_{D'}(z',x'_{k+1}) \leq
(1+2a_1)^{4a^2c_0c'M^2}e^{C+4a^2c_0M+4a^2c_0CM}\;d_{D'}(z').\eeq

The combination of  (\ref{eq(0-7-1)}) and (\ref{eq(0-7-1')}) shows
that Lemma \ref{lem13} $(3)$ holds. \epf

The following two results easily follow from the similar reasoning
as in the proof of Lemma~\ref{lem13}.

\begin{cor}\label{cor13}
For each $k\in \{1,\cdots,s\}$ and $z'\in
\beta'[x'_{1,k},x'_{1,k+1}]$,

\bee \item $d(x'_{1,k+1})\leq a_2\;d_{D'}(z')$;

\item $\lambda_{D'}(x'_{1,k+1},x'_{1,k})\leq  a_2\;d_{D'}(z')$ and

\item  $\max\{\lambda_{D'}(x'_{1,k},z'), \lambda_{D'}(x'_{1,k+1},z')\}\leq a_2\;d_{D'}(z')$. \eee
\end{cor}

\begin{cor}\label{cor13'}
For  each $z'\in \beta'[x'_{m+1},x'_{1,s+1}]$,

\bee \item $d_{D'}(w'_0)\leq a_2\;d_{D'}(z')$;

\item  $\lambda_{D'}(x'_{m+1},x'_{1,s+1})\leq
a_2\;d_{D'}(z')$ and

\item $\max\{\lambda_{D'}(x'_{m+1},z'), \lambda_{D'}(x'_{1,s+1},z')\}\leq a_2\;d_{D'}(z')$. \eee
\end{cor}
\medskip

\blem \label{lem13''} For each  $z'\in \beta'[x', w'_0]$,
$ \lambda_{D'}(x',z')\leq a_3\;d_{D'}(z'),$ and for  each $z'\in \beta'[y',w'_0]$,
$\lambda_{D'}(y',z')\leq a_3\;d_{D'}(z').$
\noindent where $a_3=a_2+a_2^2$. \elem

\bpf We only need to prove the first part of the Lemma as for the proof of the second part is similar.

If $z'\in \beta'[x', x'_{m+1}]$, then there exists some $k\in
\{1,\cdots,m\}$ such that $z'\in \beta'[x'_k, x'_{k+1}]$. If $k=1$,
then the result easily follows from Lemma \ref{lem13}. If $k>1$,
then by Lemma \ref{lem13},
\begin{eqnarray*}
\lambda_{D'}(x',z')&\leq& \lambda_{D'}(x'_1,x'_2)+\cdots + \lambda_{D'}(x'_{k-1},x'_k)+\lambda_{D'}(x'_k,z')\\
\nonumber &\leq& a_2\big(d_{D'}(x'_1)+\cdots+d_{D'}(x'_{k-1})+d_{D'}(z')\big) \\
\nonumber &\leq& (a_2+\frac{1}{2}a_2^2)d_{D'}(z').\end{eqnarray*}

Now we consider the case $z'\in \beta'[x'_{m+1}, w'_0 ]$. Then we
infer from Lemma \ref{lem13} and Corollary \ref{cor13'} that
\begin{eqnarray*}
\lambda_{D'}(x',z')&\leq& a_2\big(d_{D'}(x'_1)+d_{D'}(x'_2)+\cdots+d_{D'}(x'_{m})+d_{D'}(z')\big) \\
\nonumber &\leq& a_2\big(d_{D'}(x'_{m+1})+d_{D'}(z')\big)\\ \nonumber &\leq&
(a_2+a_2^2)d_{D'}(z').\end{eqnarray*} Hence the lemma holds.\epf

%Similarly, we have
%
%\begin{cor}\label{cor13''}
%For each $z'\in \beta'[y',w'_0]$,
%
%$$\lambda_{D'}(y',z')\leq a_3\;d_{D'}(z').$$
%where $a_3$ is the same as in Lemma~\ref{lem13''}.\end{cor}

Since $D$ is an $a$-John domain, then there exists an
$a$-cone arc $\alpha$ in $D$ joining $x$ and $y$. Let $s_0$ bisect
$\alpha$. %Then
%
%
%\begin{lem}\label{lem-4-1} Let $u\in \alpha[x, s_0]$ and $v\in \alpha[s_0, y]$. Then for each $z\in \alpha[u, s_0]$,
%$d_D(z)\geq \frac{2\ell(\alpha[u,z])+d_D(u)}{4a}$, and for each
%$z\in \alpha[s_0, v]$, $d_D(z)\geq
%\frac{2\ell(\alpha[v,z])+d_D(v)}{4a}$.
%\end{lem}

%\bpf It suffices to prove the first statement since the proof for
%the second one is similar. \epf

\begin{lem}\label{eq-8}
For each $s_1$, $s_2\in\alpha[x,s_0]$ $($or
$\alpha[x,s_0]$$)$ with $s_2\in \alpha[s_1,s_0]$, we have $$k_D(s_1,s_2)\leq  \ell_k(\alpha[s_1,s_2])\leq 2a \log\big(1+\frac{2\ell(\alpha[s_1,s_2])}{d_D(s_1)}\big)$$ and $$\ell_{k}(\alpha[s_1,s_2])\leq
4a^2c_0k_{D}(s_1,s_2)+4a^2c_0.$$ That is, for each $s_1$, $s_2\in\alpha[x,s_0]$ $($or
$\alpha[x,s_0]$$)$ \eqref{eq-3} holds.
\end{lem}

\bpf It suffices to prove the case $s_1$, $s_2\in\alpha[x,s_0]$  since the proof for the
other one is similar. For a given $s_2\in \alpha[s_1,s_0]$,
$d_D(s_2)\geq \frac{\ell(\alpha[s_1,s_2])}{a}$. If $\alpha[s_1,s_2]\subset
\mathbb{B}(s_1, \frac{d_D(s_1)}{2})$, then $d_D(z)\geq
\frac{d_D(u)}{2}$.
Otherwise, we have
$d_D(s_2)\geq\frac{d_D(s_1)}{2a}$. Hence $d_D(s_2)\geq
\frac{2\ell(\alpha[s_1,s_2])+d_D(s_1)}{4a}$, which yields that
\begin{eqnarray*}k_{D}(s_1,s_2) &\leq& \ell_k(\alpha[s_1,s_2])= \int_{\alpha[s_1,s_2]}\frac{|dz|}{d_D(z)}\\ \nonumber
&\leq& 2a\log\Big(1+\frac{2\ell(\alpha[s_1,s_2])}{d_D(s_1)}\Big)\\ \nonumber
&\leq& 4a^2\log\Big(1+\frac{d_D(s_2)}{d_D(s_1)}\Big)\\ \nonumber &\leq&
4a^2c_0k_{D}(s_1,s_2)+4a^2c_0,\end{eqnarray*} from which the proof
follows.\epf

Let $d_{D'}(v'_1)=\max\{d_{D'}(u'): u'\in\alpha'[x',s'_0]\}$ and
$d_{D'}(v'_2)=\max\{d_{D'}(u'): u'\in\alpha'[y',s'_0]\}$. Hence it follows
from Lemma \ref{lem13''} that

\blem \label{lem13-0-0} $(1)$  For each  $z'\in \alpha'[x', v'_1]$,
$\lambda_{D'}(x',z')\leq a_3\;d_{D'}(z')$
 and for each  $z'\in \alpha'[v'_1, s'_0]$, $\lambda_{D'}(s'_0,z')\leq a_3\;d_{D'}(z')$.\\
 $(2)$ For each  $z'\in \alpha'[y', v'_2]$,
$\lambda_{D'}(y',z')\leq a_3\;d_{D'}(z')$
 and for each  $z'\in \alpha'[v'_2, s'_0]$, $\lambda_{D'}(s'_0,z')\leq a_3\;d_{D'}(z')$.
\elem

%Similarly,
%
%\bcor \label{cor13-0} \ecor

\bigskip

\section{The proof of Theorem \ref{thm1}}\label{sec4}

Let $z_1$, $z_2\in D$  and $\gamma$ be a $c_0$-neargeodesic joining
$z_1$, $z_2$ in $D$. In the following, we prove that $\gamma$ is a
$b$-cone arc, that is, for each $y\in\gamma$,
\beq\label{eq(john)} \min\{\ell(\gamma[z_1, y]),\; \ell(\gamma[z_2,
y])\}\leq b\,d_D(y),\eeq where $b=4a_4c_0e^{a_4c_0}$, $a_4=a_5^{2c'M}$,
$a_5=a_6^{4a^2c_0M+C}$ and $a_6=(8a_3)^{4c'c_0M}a^2e^{2C}$.
 It is no loss of
generality to assume that $d_D(z_1)\leq d_D(z_2)$.

Let $x_0\in \gamma[z_1,z_2]$ be such that
$$d_D(x_0)=\max\limits_{z\in \gamma[z_1, z_2]}d_D(z).
$$ Then
there exists an integer $t_1 \geq 0$ such that
$$2^{t_1}\, d_D(z_1)
\leq d_D(x_0)< 2^{t_1+1}\, d_D(z_1).
$$
Let $y_0$ be the first point in $\gamma[z_1,x_0]$ from $z_1$ to
$x_0$ with
$$d_D(y_0)=2^{t_1}\, d_D(z_1).
$$
Observe that if $d_D(x_0)=d_D(z_1)$, then $y_0=z_1=x_0$.

Let $y_1=z_1$. If $z_1=y_0$, we let $y_2=x_0$. It is possible that
$y_2=y_1$. If $z_1\not= y_0$, then we let $y_2,\ldots ,y_{t_1+1}$ be
the points such that for each $i\in \{2,\ldots,t_1+1\}$, $y_i$
denotes the first point in $\gamma[z_1,x_0]$ from $y_1$ to $x_0$
satisfying
$$d_D(y_i)=2^{i-1}\, d_D(y_1).$$
Then $y_{t_1+1}=y_0$. We let $y_{t_1+2}=x_0$. It is possible that
$y_{t_1+2}=y_{t_1+1}=x_0=y_0$. This possibility occurs once
$x_0=y_0$.

Now we are going to prove for each $i\in \{1,\ldots, t_1+1\}$,
\beq\label{eq-0} k_{D}(y_i,y_{i+1})\leq a_4.\eeq

If $d_D(y_i)> \lambda_D(y_i, y_{i+1}),$ then $k_D(y_i, y_{i+1})\leq 2.$ Inequality \eqref{eq-0} obviously holds.  Hence in the following, we assume that \beq\label{eq(4-2)}d_D(y_i)\leq \lambda_D(y_i, y_{i+1}).\eeq
To prove \eqref{eq-0}, we let $\alpha_i$ be an
$a$-cone arc joining $y_i$ and $y_{i+1}$ in $D$ and let $v_i$ bisect the arclength of
$\alpha_i$. Without loss of generality, we may assume that
$d_{D'}(y'_i)\leq d_{D'}(y'_{i+1})$.

 Hence Lemma \ref{eq-8} implies
\beq\label{hl-eq(4-1-2)}k_{D}(y_i,y_{i+1})&\leq&
k_{D}(y_i,v_i)+k_{D}(y_{i+1},v_i)\\ \nonumber &\leq& 2a\Big(\log
\Big( 1+\frac{2\ell(\alpha_i[y_{i+1},v_i])} {d_D(y_{i+1})}\Big)\\
\nonumber&& + \log \Big( 1+\frac{2\ell(\alpha_i[y_i,v_i])}
{d_D(y_i)}\Big)\Big)\\ \nonumber &\leq& 4a\log \Big(
1+\frac{\ell(\alpha_i)} {d_D(y_i)}\Big). \eeq

%\blem \label{eq-0} $k_{D}(y_i,y_{i+1})\leq a_4$.\elem

Now we divide the proof of \eqref{eq-0} into two cases.

\bca \label{ca1} $\ell(\alpha_i)< a_5 \lambda_D(y_i, y_{i+1}).$\eca
Then (\ref{hl-eq(4-1-2)}) yields
\beq\label{eq(h-h-4-2')}
 \frac{\ell(\gamma[y_i,y_{i+1}])}{2d_D(y_i)}&\leq& \ell_k(\gamma[y_i,y_{i+1}])\leq c_0 k_{D}(y_i,y_{i+1})
  \leq 4ac_0\log \Big( 1+\frac{\ell(\alpha_i)} {d_D(y_i)}\Big) \\ \nonumber
  &\leq&  4ac_0\log \Big( 1+\frac{ a_5\lambda_{D}(y_i,y_{i+1})} {d_D(y_i)}\Big).\eeq A necessary condition for \eqref{eq(h-h-4-2')} is
 $$ \lambda_{D}(y_i,y_{i+1})\leq a_5^2\,d_D(y_i).$$
Hence (\ref{eq(h-h-4-2')}) implies that $k_{D}(y_i,y_{i+1})\leq
a_4$.

\bca \label{ca2} $\ell(\alpha_i)\geq a_5 \lambda_D(y_i, y_{i+1}).$\eca

We prove this case by a method of contradiction. Suppose on the contrary that
\beq\label{eq(h-4-2)}k_{D}(y_i,y_{i+1})> a_4,\eeq which implies that
\beq\label{eq(cla-0)}k_{D'}(y'_i,y'_{i+1})>\frac{k_{D}(y_i,y_{i+1})-C}{M}> 1.\eeq
Then
\begin{eqnarray*}a_4<k_{D}(y_i,y_{i+1})\leq M k_{D'}(y'_i,y'_{i+1})+C
\leq c'M\log\Big(1+\frac{\lambda_{D'}(y'_i,y'_{i+1})}{d_{D'}(y'_i)}\Big)+C,\end{eqnarray*}
and so
\beq\label{eq(h-4-1')}\lambda_{D'}(y'_i,y'_{i+1})\geq a_5d_{D'}(y'_i).\eeq
Hence
$$d_D(v_i)\geq \frac{\ell(\alpha_i)}{2a}\geq \frac{a_5}{2a}\lambda_{D}(y_i,y_{i+1})>a_6\,\lambda_{D}(y_i,y_{i+1}).$$
Then  \eqref{eq(4-2)} guarantees that there exists
$v_{i,0}\in \alpha_i[y_i,v_i]$ such that
\be\label{eq-11} d_D(v_{i,0})=a_6\,\lambda_{D}(y_i,y_{i+1}).\ee
\bcl\label{claim1}$k_{D}(y_i,v_{i,0})\leq \frac{1}{
a_5}k_{D}(y_i,y_{i+1}).$\ecl
We prove this claim also by a contradiction. Suppose that $$k_{D}(y_i,v_{i,0})> \frac{1}{
a_5}k_{D}(y_i,y_{i+1}).$$
Then Lemma \ref{eq-8} yields,
 \begin{eqnarray*}\frac{\ell(\gamma[y_i,y_{i+1}])}{2d_D(y_i)}&\leq& \ell_k(\gamma[y_i,y_{i+1}])
 \leq c_0 k_{D}(y_i,y_{i+1})\leq  a_5c_0 k_{D}(y_i,v_{i,0})\\ \nonumber&\leq& 4aa_5c_0\log\Big(1+\frac{\ell(\alpha_i[y_{i},v_{i,0}])}{d_D(y_i)}\Big)
\leq 4aa_5c_0\log\Big(1+\frac{ad(v_{i,0})}{d_D(y_i)}\Big)
\\ \nonumber
&\leq&
4a^2a_5a_6\log\Big(1+\frac{\lambda_{D}(y_i,y_{i+1})}{d_D(y_i)}\Big),
\end{eqnarray*}
whence $$ \lambda_{D}(y_i,y_{i+1})\leq a_5^2\,d_D(y_i),$$
which shows that $k_{D}(y_i,y_{i+1})\leq
a_4$ and this contradicts with $(4.6).$\qed

By \eqref{eq(4-2)} and
\eqref{eq-11}, we get
\begin{align*}
k_D(y_i, v_{i,0}) \geq \log\frac{d_D(v_{i,0})}{d_D(y_i)}\geq \log
a_6>C.\end{align*} Thus Claim \ref{claim1}, Lemma \ref{eq-8}, \eqref{eq(cla-0)} and \eqref{eq(h-4-1')}
imply that
\begin{eqnarray*} \log \Big(1+\frac{\lambda_{D'}(y'_i,v'_{i,0})}{d_{D'}(y'_i)}\Big)&\leq& k_{D'}(y'_i,v'_{i,0})
\leq M k_{D}(y_i,v_{i,0})+C \\ &<& 2Mk_{D}(y_i,v_{i,0})
\leq\frac{2M}{
a_5}k_{D}(y_i,y_{i+1})\\  &\leq&\frac{2M^2}{ a_5}k_{D'}(y'_i,y'_{i+1})+\frac{2CM}{ a_5}\\
&\leq&\frac{4c'M^2}{ a_5}\log\Big(1+\frac{\lambda_{D'}(y'_i,y'_{i+1})}{d_{D'}(y'_i)}\Big)\\
 &\leq&\log
\Big(1+\frac{\lambda_{D'}(y'_i,y'_{i+1})}{a_5d_{D'}(y'_i)}\Big).\end{eqnarray*}
Hence
\beq\label{eq(hl-41-5)}\lambda_{D'}(y'_i,v'_{i,0})<
\frac{1}{a_5}\lambda_{D'}(y'_i,y'_{i+1}),\eeq which, together with
(\ref{eq(h-4-1')}), gives
\beq\label{eq--2} d_{D'}(v'_{i,0})\leq
\lambda_{D'}(y'_i,v'_{i,0})+d_{D'}(y'_i)\leq
\frac{2}{a_5}\lambda_{D'}(y'_i,y'_{i+1}).\eeq
\bcl\label{eq--6} $\lambda_{D'}(y'_i,v'_i)<
\frac{\lambda_{D'}(y'_i,y'_{i+1})}{2}$. \ecl

Suppose on the contrary that $$\lambda_{D'}(y'_i,v'_i)\geq
\frac{\lambda_{D'}(y'_i,y'_{i+1})}{2}.$$ Let $u'_{0,i}\in\gamma'[y'_{i},
y'_{i+1}]$ be a point satisfying $$d_{D'}(u'_{0,i})=\max\{d_{D'}(w'):w'\in\gamma'[y'_{i},
y'_{i+1}]\}.$$ Obviously, $$\max\{\lambda_{D'}(y'_{i+1},u'_{0,i}),
\lambda_{D'}(u'_{0,i},y'_i)\}\geq \frac{\lambda_{D'}(y'_i,y'_{i+1})}{2}.$$
Then we know from Lemma \ref{lem13''}
that
\be\label{e---1} d_{D'}(u'_{0,i})\geq
\frac{\lambda_{D'}(y'_i,y'_{i+1})}{2a_3}.\ee Hence by Lemma
\ref{lem13''} and (\ref{eq(h-4-1')}), there must exist some point
$y'_{0,i}\in \gamma'[y'_i,u'_{0,i}]$ satisfying
\beq\label{eq(W-l-6-1)}d_{D'}(y'_{0,i})=\frac{\lambda_{D'}(y'_i,y'_{i+1})}{2a_3}
\,\; \mbox{and}\, \;\lambda_{D'}(y'_i,y'_{0,i})\leq
a_3\,d_{D'}(y'_{0,i}).\eeq

Let $v'_0\in\alpha'_i[y'_{i}, v'_{i}]$ satisfy $d_{D'}(v'_0)=\max\{d_{D'}(u'):u'\in\alpha'_i[y'_{i}, v'_{i}]\}$.
Then Lemma
\ref{lem13-0-0} shows that for each $z'\in \alpha'_i[ v'_{0}, v'_i]$, \beq\label{cla-3}\lambda_{D'}(v'_i,z')\leq a_3 d_{D'}(z').\eeq
By \eqref{eq(hl-41-5)} and \eqref{eq--2} we have
\begin{eqnarray*}\lambda_{D'}(v'_i,v'_{i,0})&\geq&
\lambda_{D'}(v'_i,y'_i)-\lambda_{D'}(v'_{i,0},y'_i)\\ \nonumber
&\geq& (\frac{1}{2}-\frac{1}{a_5})\lambda_{D'}(y'_i,y'_{i+1})\\ \nonumber
&\geq& (a_5-\frac{1}{2})d_{D'}(v'_{i,0}),
\end{eqnarray*}
which together with \eqref{cla-3} show that $v'_0\in \alpha'_i[ v'_{i,0}, v'_i]$.

Obviously, $\max\{\lambda_{D'}(v'_{i},v'_0),
\lambda_{D'}(v'_0,y'_i)\}\geq
\frac{\lambda_{D'}(y'_i,y'_{i+1})}{4}$. We know from Lemma
\ref{lem13-0-0} that $d_{D'}(v'_0)\geq
\frac{\lambda_{D'}(y'_i,y'_{i+1})}{4a_3}$.  By (\ref{eq--2}) and
Lemma \ref{lem13-0-0}, we see that there exists some point $u'_0\in
\alpha'_i[v'_{i,0},v'_{i}]$ such that
\beq\label{eq(W-l-6-2)}
d_{D'}(u'_0)=\frac{\lambda_{D'}(y'_i,y'_{i+1})}{4a_3}\,\; \mbox{and}\,
\;\lambda_{D'}(y'_i,u'_0)\leq a_3\,d_{D'}(u'_0).\eeq Hence
(\ref{eq(W-l-6-1)}) shows that
\begin{eqnarray*}
\log \frac{d_D(u_0)}{d_D(y_{0,i})}&\leq& k_{D}(y_{0,i},u_0)
\leq M k_{D'}(y'_{0,i},u'_0)+C \\ \nonumber
&\leq&Mc'\log\Big(1+\frac{\lambda_{D'}(u'_0,y'_{0,i})}{\min\{d_{D'}(u'_0),
d_{D'}(y'_{0,i})\}}\Big)+C\\ \nonumber
&\leq&Mc'\log\Big(1+\frac{\lambda_{D'}(u'_0,y'_i)+\lambda_{D'}(y'_i,y'_{0,i})}{\min\{d_{D'}(u'_0),
d_{D'}(y'_{0,i})\}}\Big)+C
\\ \nonumber
&\leq&
Mc'\log\Big(1+\frac{a_3d(u'_0)+\lambda_{D'}(y'_i,y'_{0,i})}{\min\{d_{D'}(u'_0),
d_{D'}(y'_{0,i})\}}\Big)+C\\ \nonumber &<&Mc'\log (1+3a_3)+C,
\end{eqnarray*}
which yields that
\beq\label{eq(W-l-6-4)}d_D(u_0)\leq (1+3a_3)^{Mc'}e^Cd_D(y_{0,i}).\eeq
Lemma \ref{eq-8}, \eqref{eq--2} and \eqref{eq(W-l-6-2)} imply
that
\begin{eqnarray*} 4{a}^2M\log\Big(1+\frac{d_D(u_0)}{d_D(v_{i,0})}\Big)+C
 &\geq& M\ell_k(\alpha_i[v_{i,0},u_0])+C
 \geq  M k_{D}(v_{i,0},u_0)+C\\  &\geq&k_{D'}(v'_{i,0},u'_0)
\geq\log\frac{d_{D'}(u'_0)}{d_{D'}(v'_{i,0})}
\geq\log\frac{a_5}{8a_3},\end{eqnarray*} whence $d_D(u_0)\geq
a_6d_D(v_{i,0})$. So we infer from \eqref{eq(4-2)} and
\eqref{eq-11} that
\begin{align*} d_D(u_0)\geq  a_6d_D(v_{i,0})=a_6^2\lambda_D(y_i, y_{i+1})
\geq \frac{ a_6^2}{2}d_D(y_{i+1})\geq \frac{
a_6^2}{2}d_D(y_{0,i}),\end{align*} which contradicts with
(\ref{eq(W-l-6-4)}). Hence Claim \ref{eq--6} holds.\qed
\medskip

It is obvious from Claim \ref{eq--6} that
$\lambda_{D'}(y'_{i+1},v'_i)>
\frac{\lambda_{D'}(y'_i,y'_{i+1})}{2}$. Let $q'_0\in
\alpha'_i[y'_i,v'_i]$ with
\beq\label{112}\frac{\lambda_{D'}(y'_i,v'_i)}{2a_3}\geq\lambda_{D'}(q'_0,v_i')\geq\frac{\lambda_{D'}(y'_i,v'_i)}{4a_3},\eeq
and $u'_1\in\alpha'_i[y'_{i+1},v'_i]$ with
\beq\label{122}\frac{\lambda_{D'}(y'_i,v'_i)}{2a_3}\geq\lambda_{D'}(u'_1,v_i')\geq\frac{\lambda_{D'}(y'_i,v'_i)}{4a_3}.\eeq
By Lemma \ref{lem13-0-0}, we  get
\beq\label{eq(W-l-6'-0)} d_{D'}(q'_0)\geq
\frac{\lambda_{D'}(y'_i,v'_i)}{4a^2_3} \,\; \mbox{and}\;\,
\label{eq(W-l-6'-1)}d_{D'}(u'_1)\geq
\frac{\lambda_{D'}(y'_i,v'_i)}{4a_3^2}.\eeq Hence we have
\beq \label{eq--7} \Big|\log\frac{d_D(u_1)}{d_D(q_0)}\Big| &\leq&k_{D}(u_1,q_0)\\ \nonumber&\leq& M k_{D'}(u'_1, q'_0)+C\\
\nonumber &\leq&Mc'\log\Big(1+\frac{\lambda_{D'}(u'_1,
q'_0)}{\min\{d_{D'}(q'_0), d_{D'}(u'_1)\}}\Big)+C\\
\nonumber &\leq& Mc'\log\Big(1+\frac{\lambda_{D'}(u'_1,v'_i
)+\lambda_{D'}(v'_i,q'_0)}{\min\{d_{D'}(q'_0), d_{D'}(u'_1)\}}\Big)+C\\
\nonumber &\leq& Mc'\log (1+4a_3 )+C,\eeq which implies that
\beq\label{eq(W-l-6'-2)}\frac{d_D(u_1)}{(1+4a_3 )^{Mc'}e^C}\leq
d_D(q_0)\leq (1+4a_3 )^{Mc'}e^Cd_D(u_1).\eeq

\bcl\label{eq(W-l-6'-3)} $d_D(q_0)\geq a_5d_D(v_{i,0})$.\ecl
\noindent Otherwise, Lemma \ref{eq-8}, (\ref{eq-11}),
\eqref{eq--7} and (\ref{eq(W-l-6'-2)})
 show that
\beq\label{eq(W-l-6'-2')}\frac{\ell(\gamma[y_i,y_{i+1}])}{2d_D(y_i)}&\leq&
\ell_k(\gamma[y_i,y_{i+1}])
\leq
c_0 k_{D}(y_i, y_{i+1})\\
\nonumber&\leq& c_0 (k_{D}(y_i, q_0)+k_{D}(q_0, u_1)+k_{D}(u_1, y_{i+1}))\\
\nonumber &\leq& 4a^2c_0\log\Big(1+\frac{d_D(
q_0)}{d_D(y_i)}\Big)+Mc'c_0\log\Big(1+2a_3\Big)\\
\nonumber &&+Cc_0+4a^2c_0\log\Big(1+\frac{d_D(
u_1)}{d_D(y_{i+1})}\Big)\\
\nonumber &\leq& 9a^2a_5c_0\log\Big(1+\frac{\lambda_D(y_i,
y_{i+1})}{d_D(y_i)}\Big).\eeq A necessary condition for
\eqref{eq(W-l-6'-2')} is $\lambda_D(y_i, y_{i+1})\leq a_5^2d_D(y_i)$.
Hence by (\ref{eq(W-l-6'-2')}), we know that
$$k_{D}(y_i, y_{i+1})\leq 9a^2a_5\log(1+a_5^2),$$ which contradicts with (\ref{eq(h-4-2)}).
We complete the proof of Claim \ref{eq(W-l-6'-3)}.\qed
\medskip

By \eqref{eq(h-4-1')} and \eqref{e---1}
$$\lambda_{D'}(u'_{0,i},y'_i)\geq d_{D'}(u'_{0,i})-d_{D'}(y'_i)\geq \frac{1}{3a_3}\lambda_{D'}(y'_{i+1},y'_i).$$
Then Claim \ref{eq--6} guarantees that there exists $y'_0\in
\gamma'[y'_i,u'_{0,i}]$ such that
\beq\label{132}\frac{\lambda_{D'}(y'_i,v'_i)}{2a_3}\geq\lambda_{D'}(y'_0,y'_i)\geq\frac{\lambda_{D'}(y'_i,v'_i)}{3a_3}.\eeq
 Hence Lemma \ref{lem13''}
implies that
$$\frac{\lambda_{D'}(y'_i,v'_i)}{2a_3}=\lambda_{D'}(y'_0,y'_i)\leq a_3 d_{D'}(y'_0).$$ Hence \eqref{112}, \eqref{122}, \eqref{eq(W-l-6'-0)} and \eqref{132} give
\begin{eqnarray*}\log \frac{d_D(q_0)}{d_D(y_{0})}&\leq& k_{D}(q_0, y_{0})
\leq M k_{D'}(q'_0,y'_0)+C\\
\nonumber &\leq&
Mc'\log\Big(1+\frac{\lambda_{D'}(y'_0,q'_0)}{\min\{d_{D'}(q'_0),d_{D'}(y'_0)\}}\Big)+C
\\
\nonumber &\leq&
Mc'\log\Big(1+\frac{\lambda_{D'}(y'_i,v'_i)+\lambda_{D'}(v'_i,q'_0)+\lambda_{D'}(y'_i,y'_0)}{\min\{d_{D'}(q'_0),d_{D'}(y'_0)\}}\Big)+C
\\
\nonumber &\leq& Mc'\log(1+4a_3+4a_3^2)+C.\end{eqnarray*} We infer
from  \eqref{eq(4-2)} and \eqref{eq-11} that
\begin{eqnarray*}d_D(q_0)&\leq&
(1+4a_3+4a_3^2)^{Mc'}e^C d_D(y_0)\\ \nonumber&\leq&
2(1+4a_3+4a_3^2)^{Mc'}e^C d_D(y_i)
\\ \nonumber&\leq& 2(1+4a_3+4a_3^2)^{Mc'}e^C \lambda_{D}(y_i,y_{i+1})\\
\nonumber&=& \frac{ 2(1+4a_3+4a_3^2)^{Mc'}e^C}{a_6}
d_D(v_{i,0}),\end{eqnarray*} which contradicts with Claim
\ref{eq(W-l-6'-3)}. We complete the proof of \eqref{eq-0}.

\medskip

Then by  \eqref{eq-0} we have for all $i\in\{1,\cdots, t_1+1\}$,
\beq\label{eq(h-56)} \frac{\ell(\gamma[y_i,y_{i+1}])}{2d_D(y_i)}\leq
\ell_k(\gamma[y_i,y_{i+1}])\leq c_0  k_{D}(y_i,y_{i+1}) \leq a_4c_0,\eeq
which implies that
\be\label{h-57} \ell(\gamma[y_i,y_{i+1}])\leq
 2a_4c_0\,d_D(y_i).\ee

Further, for each $y\in \gamma[y_i,y_{i+1}]$, it follows from
\eqref{eq(h-56)} that
\beq\label{eq(h-58)} \log \frac{d_D(y_i)}{d_D(y)}\leq k_D(y,y_i)\leq
\,c_0k_D(y_i,y_{i+1})\leq a_4c_0,\eeq whence
$$d_D(y_i)\leq e^{ a_4c_0}d_D(y).$$
For each $y\in \gamma[y_1,x_{0}]$, there is some $i\in
\{1,\cdots,t_1+1\}$ such that $y\in \gamma[y_i,y_{i+1}]$. It follows
from (\ref{h-57}) and (\ref{eq(h-58)}) that
\beq\label{eq(h-59)} \ell(\gamma[z_1,y])&=&
\ell(\gamma[y_1,y_2])+\ell(\gamma[y_2,y_3])+\cdots+\ell(\gamma[y_i,y])
\\
\nonumber &\leq& 2a_4c_0(d_D(y_1)+d_D(y_2)+\cdots+d_D(y_i))\\ \nonumber
&\leq& 4a_4c_0\,d_D(y_i)\\ \nonumber &\leq& 4a_4c_0e^{a_4c_0}\,d_D(y).\eeq

By replacing $\gamma[z_1, x_0]$ by $\gamma[z_2, x_0]$ and repeating
the procedure as above, we also get that
\beq\label{eq(h-59')} \ell(\gamma[z_{2},y])\leq
4a_4c_0e^{ a_4c_0}\,d_D(y).\eeq The combination of \eqref{eq(h-59)} and
\eqref{eq(h-59')} concludes the proof of Theorem \ref{thm1}.\qed

\bigskip
{\sc Acknowledgements.}
This work was partially supported by a grant from Simons Foundation and Talent Introduction Foundation of Central South University of Forestry and Technology (No. 2013RJ005). The author would like to thank Professor Manzi Huang for several comments on this manuscripts and thank the referee who have made valuable comments on this manuscripts. The revision of this work was completed during the visit of the author to Indian Statistical Institute (ISI) Chennai Centre, India. She thanks ISI Chennai Centre for the hospitality and the other supports.

\end{document}